\documentclass[12pt]{article}

\usepackage{xspace}
\usepackage{enumerate}
\usepackage{latexsym,epsfig}
\usepackage{amssymb}
\usepackage{amsthm}
\usepackage{amsmath}

\newcommand{\la}{\lambda}

\newcommand{\C}{\cot_{p}}
\newcommand{\T}{\tan_{p}}
\newcommand{\Si}{\sin_{p}}
\newcommand{\Sh}{\sinh_{p}}

\title{The dual eigenvalue problems for $p$-Laplacian}

\author{Y.H. Cheng$^1$, Wei-Cheng
Lian$^2$ and Wei-Chuan Wang$^3$}

\date{\today}

\begin{document}

\maketitle

\begin{abstract}
In this paper, we find the minimizer of the eigenvalue gap for the
single-well potential problem and the eigenvalue ratio for the
single-barrier density problem and symmetric single-well
(single-barrier)density problem for $p$-Laplacian. This extends the
results of the classical Sturm-Liouville problem.
\end{abstract}


\footnote{$p$-Laplacian; eigenvalue gap; eigenvalue ratio}

\footnote{$^1$Department of Mathematics, National Tsing Hua University, Hsinchu, Taiwan 300, R.O.C. Email:
jengyh@math.nsysu.edu.tw}

\footnote{$^2$Department of Information Management, National
 Kaohsiung Marine Univeristy, Kaohsiung, Taiwan, R.O.C. Email:
 wclian@mail.nkmu.edu.tw}

\footnote{$^3$Department of Applied Mathematics, National Sun Yat-sen
University, Kaohsiung, Taiwan 804, R.O.C. Email:
wangwc@math.nsysu.edu.tw}

\section{Introduction}
For $q,\ \rho\in L^{1}, \ \rho >0$ a.e., and $p>1$,  consider the eigenvalue problem for $p$-Laplacian
\begin{equation}(y'^{(p-1)})'=-(p-1)(\la\rho(x) -q(x))y^{(p-1)}\ ,
\label{eq1.1}
\end{equation}
with the Dirichlet boundary conditions \begin{equation}y(0)=y(\hat{\pi})=0\ .\label{eq1.2}
\end{equation}
For $p=2$, (\ref{eq1.1}) is reduced to Schr\"odinger equation  $y''=-(\lambda -q(x)) y$ when $\rho\equiv 1$, while  (\ref{eq1.1}) is reduced to the string equation  $y''=-\lambda\rho(x) y$ when $q\equiv 0$.

Denote by $\Si (x)$ the solution of
\begin{equation}
\label{eq3.1}
\left\{\begin{array}{l} (y'^{(p-1)})'=-(p-1)y^{(p-1)}\ ,\\ y(0)=0\ ,\ y'(0)=1\ . \end{array}\right.
\end{equation}
Then we have
\begin{equation}
\label{eq2.3}
|\Si (x)|^{p}+|\Si' (x)|^{p}=1\ .
\end{equation}
Here, $\Si (x)$ is called a general sine function. In \cite{E79},
Elbert discussed the analogies between $\Si (x)$ and $\sin x$.
For example, he showed that $w=w(x)=\Si (x)$ is the inverse
function of the below integral
$$x=\int_{0}^{w}\frac{dt}{(1-t^{p})^{\frac{1}{p}}}\ ,\
\mbox{for}\ 0\leq w\leq 1\ ,$$ and $\Si(x) = 1$ at
$x=\frac{\hat{\pi}}{2}\equiv
\int_{0}^{1}\frac{dt}{(1-t^{p})^{\frac{1}{p}}}=\frac{\pi}{p
\sin(\pi/p)}$. Furthermore, defining
$$\Si(x)=\left\{\begin{array}{ll}\Si(\hat{\pi}-x)\ ,& \mbox{if}\ \frac{\hat{\pi}}{2}\leq x\leq \hat{\pi}\ ,\\
-\Si(x-\hat{\pi})\ ,& \mbox{if}\ \hat{\pi}\leq x\leq 2\hat{\pi}\ ,\\
\Si(x-2n\hat{\pi})\ ,& \mbox{for}\ n=\pm 1, \pm 2, \cdots\ ,\\
 \end{array} \right.$$
he obtained a sine-like function.
Note that $\hat{\pi}$ is the first zero of $\Si (x)$.

Recently,  there have been a number of studies on the optimal
estimates of eigenvalues, eigenvalue gaps and eigenvalue ratios for
eigenvalue problem $-y''+q(x)y=\lambda\rho(x) y$
\cite{K76,AS92,L94,CS94,H02,H07}. It was proved that, for
Schr\"odinger equation $-y''+q(x)y=\lambda y$, the constant
potential function gives the minimum Dirichlet eigenvalue gaps
$\lambda_{2}-\lambda_{1}$ when the potential function $q$ is assumed
to be convex \cite{L94}, symmetric single-well \cite{AB89} or
single-well \cite{H02}, while under some additional conditions, the
symmetric 1-step function is the potential function in $E[h,H,M]$
giving the minimal Dirichlet eigenvalue gap \cite{CS94}. On the
other hand, it is known that, for the string equation $-y''=\la
\rho(x) y$, the constant density function gives the minimum
Dirichlet eigenvalue ratio $\frac{\lambda_{2}}{\lambda_{1}}$ when
the density function $\rho$ is assumed to be concave, symmetric
single-barrier \cite{H99} or single-barrier \cite{H02}, while the
symmetric 1-step function is the density in $E[h,H,M]$ giving the
minimum Dirichlet eigenvalue ratio \cite{K76}, see also \cite{MW76}.
These results are called "duality results". In particular, Ashbaugh
and Benguria in 1989  found the optimal bound of the eigenvalue
ratio $\la_{n}/\la_{1}$ for Schr\"odinger equation with nonnegative
potentials \cite{AB892}, and this result was extended by Huang and
Law for general Sturm-Liouville problems \cite{HL96}. It shall be
mentioned that Huang in 2007 discuss the eigenvalue gap for
vibrating string with symmetric single-well densities \cite{H07}.
Here, the function $V$ is called a single-well function with the
transition point $a$ if $V(x)$ is decreasing in $[0,a]$ and
increasing in $[a,\pi]$ while $V$ is called a single-barrier function if $-V$ is a single-well function.

In this paper, we will generalize the results of the Dirichlet
eigenvalue gap for Schr\"odinger equation and eigenvalue ratio for
string equation in \cite{H99,H02} to  $p$-Laplacian. We
obtain the following results.

\newtheorem{th1.1}{Theorem}[section]
\begin{th1.1}
\label{th1.1}

Consider the eigenvalue problem for $p$-Laplacian
(\ref{eq1.1})-(\ref{eq1.2}) with $\rho\equiv 1$.

If $q$ is single-well with a transition point at $\frac{\hat{\pi}}{2}$, then $$\la_{2}-\la_{1}\geq 2^{p}-1\ .$$
The equality holds if and only if $q$ is constant. Furthermore, if the  transition point $a\neq \frac{\hat{\pi}}{2}$, then there is a single-well potential such that $\la_{2}-\la_{1}<2^{p}-1\ .$
\end{th1.1}

\newtheorem{th1.2}[th1.1]{Theorem}
\begin{th1.2}
\label{th1.2}

Consider the eigenvalue problem for $p$-Laplacian
(\ref{eq1.1})-(\ref{eq1.2}) with $q\equiv 0$.

\begin{enumerate}
\item[(a)]
If $\rho$ is single-barrier density with a transition point at
$\frac{\hat{\pi}}{2}$, then $$\frac{\mu_{2}}{\mu_{1}}\geq 2^{p}\ .$$
The equality holds if and only if $\rho$ is constant. Furthermore, if the  transition point $a\neq \frac{\hat{\pi}}{2}$, then there is a single-barrier density such that $\frac{\mu_{2}}{\mu_{1}}< 2^{p}$.

\item[(b)]
If $\rho$ is a symmetric single-well density with a transition point
at $\frac{\hat{\pi}}{2}$, then
$$\frac{\mu_2}{\mu_1}\le 2^p.$$ The equality holds
if and only if $\rho(x)$ is a constant a.e..
\end{enumerate}
\end{th1.2}

\section{Preliminaries}
As in Binding and Drabek \cite{BD2003}, the eigenvalues $\lambda_k$,
form a strictly increasing sequence as
\begin{equation}
\lambda_1[\rho,\ q]<\lambda_2[\rho,\ q]<\lambda_3[\rho,\ q]<\cdots,
\label{eq1.3}
\end{equation}
and accumulating at $\infty$. The $n$-th eigenfunction $y_n$ has
$n-1$ zeros in $(0,\hat{\pi})$.

Let $y_{n}(x)=y(x,\la_{n})$ be the $n$-th normalized eigenfunction of
(\ref{eq1.1})-(\ref{eq1.2}) satisfying
$\int_{0}^{\hat{\pi}}\rho(x)|y(x)|^{p}dx=1$. We may assume
$y_{n}(x)>0$ initially and let $x_{0}$ be the zero of $y_{2}(x)$.  In order to compare the behaviors of
$y_1$ and $y_2$, we introduce a Pr\"ufer-type substitution. Let
\begin{equation*}
 y_n(x)=r(x) \Si(\phi _n(x))\ ,\quad y'_n(x) =
 r(x)\Si'(\phi _n(x))\ .
\end{equation*}

Denote by  $\T (x) = \frac{\Si(x)}{\Si'(x)}$ and $\C (x)=\frac{\Si'(x)}{\Si(x)}$  the generalized tangent and cotangent functions respectively. Since $$\C'(x) = \frac{d}{dx} \frac{\Si'(x)}{\Si(x)}=-\left|\frac{\Si(x)}{\Si'(x)}\right|^{p-2}-|\C(x)|^{2}=-(1+|\T(x)|^{p})|\C(x)|^{2}\ ,$$ the function $\C (x)$
is strictly decreasing on $(0,\hat{\pi})$. This implies
$$(\frac{y_2}{y_1})'=\frac{y_1y_2'-y_2y_1'}{y_1^2}=\frac{y_1y_2}{y_1^2}[\frac{y_2'}{y_2}-\frac{y_1'}{y_1}]=\frac{y_2}{y_1}
[\C(\phi _2(x))-\C(\phi _1(x))]\ .$$After the Pr\"ufer substitution, we
obtain
$$\phi _n'=|\Si'(\phi _n)|^p+(\lambda _n\rho(x)-q(x)) |\Si(\phi _n)|^p.$$ By Comparison theorem \cite{br89}, we
have $\phi _2(x)>\phi _1(x)$ on $(0,x_{0})$ and, hence,
$(\frac{y_2}{y_1})'< 0$ on $(0,x_{0})$. This implies
$\frac{y_2}{y_1}$ is strictly decreasing on $(0,x_{0})$.
Furthermore, $y_{1}$ and $y_{2}$ has at most one intersection point
in $(0,x_{0})$. Similarly, $y_{1}$ and $-y_{2}$ has at most one
intersection point in $(x_{0},\hat{\pi})$. Hence we have the
following lemma.

\newtheorem{th2.2}{Lemma}[section]
\begin{th2.2}
\label{th2.2}
Consider the eigenvalue problem for $p$-Laplacian (\ref{eq1.1})-(\ref{eq1.2}).  Then $|y_{1}(x)|=|y_{2}(x)|$ have at most two intersection points on $(0,\hat{\pi})$.
\end{th2.2}

Let $\rho(x,t)$ and $q(x,t)$ be  one-parameter family of piecewise
continuous functions such that $\frac{\partial}{\partial t}\rho$ and
$\frac{\partial}{\partial t}q$ exist. Denote by
$\{(\la_{n}(t),y_{n}(x,t))\}_{n\geq 1}$  the $n$-th normalized
eigenpair. The following lemma is an extension for the case $p=2$ in
\cite{K76} (see also \cite{L94,H99}). The proof will be given in
appendix.

\newtheorem{th2.1}[th2.2]{Lemma}
\begin{th2.1}
\label{th2.1}
\begin{equation}\frac{d}{dt}\la_{n}(t)=\int_{0}^{\hat{\pi}}\frac{\partial}{\partial t}q(x,t)|y_{n}(x,t)|^{p}dx- \la_{n}\int_{0}^{\hat{\pi}}\frac{\partial}{\partial t}\rho(x,t)|y_{n}(x,t)|^{p}dx\ .
\label{eq2.1}
\end{equation}

\end{th2.1}

Following from Lemma \ref{th2.1}, we have
\begin{enumerate}
  \item[1.] If $\rho\equiv 1$, we have
$$\frac{d}{dt}\left(\lambda_n(t)-\lambda_m(t)\right)=\int_0^{\hat{\pi}}\frac{\partial
q}{\partial
t}(x,t)\left(|y_n(x,t)|^p-|y_m(x,t)|^p\right)dx\ ;$$

  \item[2.]
If $q\equiv 0$, we have
  $$
  \frac{d}{dt}\left(\frac{\lambda_n(t)}{\lambda_m(t)}\right)=\frac{\lambda_n(t)}{\lambda_m(t)}\int_0^{\hat{\pi}}\frac{\partial
\rho}{\partial
t}(x,t)\left(|y_m(x,t)|^p-|y_n(x,t)|^p\right)dx\ .$$
\end{enumerate}

Next, Lemma \ref{th2.3} will be used to proof the eigenvalue gap (Theorem \ref{th1.1}) while Lemma \ref{th2.4} will be used to proof the eigenvalue ratio (Theorem \ref{th1.2}).

\newtheorem{th2.3}[th2.2]{Lemma}
\begin{th2.3}
\label{th2.3}

Denote  $f(t)=t^{\frac{1}{p}}\C (t^{\frac{1}{p}}\frac{\hat{\pi}}{2})$. Let $t_{n}$ be the $n$-th solution of $f(t)=-f(t-m)$ where $m>0$. Then $$t_{2}-t_{1}\geq 2^{p}-1\ .$$

\end{th2.3}

\begin{proof}

Note that, according graph analysis, $t_{1}\in (1,\min\{1+m,2^{p}\})$ for $m>0$. For $m\geq 3^{p}-1$, we have $t_{2}\geq 3^{p}$ and hence
$$t_{2}-t_{1}\geq 3^{p}-2^{p}>2^{p}-1\ .$$
So we only need to consider $0<m<3^{p}-1$. In this case, $t_{2}\in (2^{p}, \min\{2^{p}+m,3^{p}\})$.


\begin{enumerate}
\item[1.] Assume $t\geq 0$.
By the definition, we have  $f(t)=t^{\frac{1}{p}}\C (t^{\frac{1}{p}}\frac{\hat{\pi}}{2})$,
\begin{eqnarray*}
f'(t)&=&\frac{1}{p}t^{\frac{1-p}{p}}\C(t^{\frac{1}{p}}\frac{\hat{\pi}}{2})-t^{\frac{1}{p}}(1+|\T(t^{\frac{1}{p}}\frac{\hat{\pi}}{2})|^{p})\C^{2}(t^{\frac{1}{p}}\frac{\hat{\pi}}{2})\cdot\frac{1}{p}t^{\frac{1-p}{p}}\frac{\hat{\pi}}{2}\ \\
&=&\frac{1}{pt}f(t)-\frac{\hat{\pi}}{2pt}(1+|\T(t^{\frac{1}{p}}\frac{\hat{\pi}}{2})|^{p})|f(t)|^{2}\ \\
&=&\frac{t^{\frac{1-p}{p}}}{2p|\Si(t^{\frac{1}{p}}\frac{\hat{\pi}}{2})|^{2}}(2\Si(t^{\frac{1}{p}}\frac{\hat{\pi}}{2})\Si'(t^{\frac{1}{p}}\frac{\hat{\pi}}{2})-t^{\frac{1}{p}}\hat{\pi}|\Si'(t^{\frac{1}{p}}\frac{\hat{\pi}}{2})|^{2-p})\ .
\end{eqnarray*}

If $\Si'(t^{\frac{1}{p}}\frac{\hat{\pi}}{2}) >0$, in this case $t^{\frac{1}{p}}\in (0,1)$ and $(4n-1,4n+1)$ for $n\geq 1$, then \begin{eqnarray*}
2\Si(t^{\frac{1}{p}}\frac{\hat{\pi}}{2})\Si'(t^{\frac{1}{p}}\frac{\hat{\pi}}{2})-t^{\frac{1}{p}}\hat{\pi}|\Si'(t^{\frac{1}{p}}\frac{\hat{\pi}}{2})|^{2-p}&=&\Si'(t^{\frac{1}{p}}\frac{\hat{\pi}}{2}) (2\Si(t^{\frac{1}{p}}\frac{\hat{\pi}}{2})-t^{\frac{1}{p}}\hat{\pi}|\Si'(t^{\frac{1}{p}}\frac{\hat{\pi}}{2})|^{1-p})\ ,\\
&\leq&\Si'(t^{\frac{1}{p}}\frac{\hat{\pi}}{2}) (2\Si(t^{\frac{1}{p}}\frac{\hat{\pi}}{2})-t^{\frac{1}{p}}\hat{\pi})\ ,\\
&\equiv& \Si'(t^{\frac{1}{p}}\frac{\hat{\pi}}{2}) g(t)\ .
\end{eqnarray*}
Since $g(0)=0, g((4n-1)^{p})$ and $g'(t)=\frac{t^{\frac{1-p}{p}}\hat{\pi}}{p}(\Si'(t^{\frac{1}{p}}\frac{\hat{\pi}}{2})-1)<0$ for  $t^{\frac{1}{p}}\in (0,1)$ and $(4n-1,4n+1), n\geq 1$, we have $g(t)<0$ for  $t^{\frac{1}{p}}\in (0,1)$ and $(4n-1,4n+1), n\geq 1$ and hence $f'(t)<0$ for $t^{\frac{1}{p}}\in (0,1)$ and $(4n-1,4n+1), n\geq 1$,.

Similarly, if  $\Si'(t^{\frac{1}{p}}\frac{\hat{\pi}}{2}) <0$,  in this case $t^{\frac{1}{p}}\in  (4n-3,4n-1)$ for $n\geq 1$, then \begin{eqnarray*}
2\Si(t^{\frac{1}{p}}\frac{\hat{\pi}}{2})\Si'(t^{\frac{1}{p}}\frac{\hat{\pi}}{2})-t^{\frac{1}{p}}\hat{\pi}|\Si'(t^{\frac{1}{p}}\frac{\hat{\pi}}{2})|^{2-p}&=&\Si'(t^{\frac{1}{p}}\frac{\hat{\pi}}{2}) (2\Si(t^{\frac{1}{p}}\frac{\hat{\pi}}{2})+t^{\frac{1}{p}}\hat{\pi}|\Si'(t^{\frac{1}{p}}\frac{\hat{\pi}}{2})|^{1-p})\ ,\\
&\leq&\Si'(t^{\frac{1}{p}}\frac{\hat{\pi}}{2}) (2\Si(t^{\frac{1}{p}}\frac{\hat{\pi}}{2})+t^{\frac{1}{p}}\hat{\pi})\ ,\\
&\equiv& \Si'(t^{\frac{1}{p}}\frac{\hat{\pi}}{2}) h(t)\ .
\end{eqnarray*}
Since $h(0)=0, h((4n-3)^{p})>0$ and $h'(t)=\frac{t^{\frac{1-p}{p}}\hat{\pi}}{p}(\Si'(t^{\frac{1}{p}}\frac{\hat{\pi}}{2})+1)>0$, we have $h(t)>0$ for $t^{\frac{1}{p}}\in  (4n-3,4n-1), n\geq 1$ and hence $f'(t)<0$ for $t^{\frac{1}{p}}\in  (4n-3,4n-1), n\geq 1$.

\item[2.] Assume $t< 0$.
Define   by $w=w(x)=\Sh (x)$  the inverse function of the integral $x=\int_{0}^{w}\frac{dt}{(1+t^{p})^{\frac{1}{p}}}$.
We call $\Sh (x)$  the generalized hyperbolic sine function.
 It is easy to show that $\Sh (x)=(-1)^{-\frac{1}{p}}\Si ((-1)^{\frac{1}{p}}x)$ and $\Sh' (x)=\Si' ((-1)^{\frac{1}{p}}x)$ where $(-1)^{\frac{1}{p}}=e^{\pi i/p}$.
 Furthermore,
\begin{equation}
\label{eq2.5}
\Sh'^{p}(x)-\Sh^{p}(x)=1\ ,
\end{equation}
 and then $\Sh''(x)=\frac{\Sh^{p-1}(x)}{\Sh'^{p-2}(x)}$.

  Let $\hat{t}=-t$. Since
$$f(t)=t^{\frac{1}{p}}\C (t^{\frac{1}{p}}\frac{\hat{\pi}}{2})=(-1)^{\frac{1}{p}}\hat{t}^{\frac{1}{p}}\frac{\Si'((-1)^{\frac{1}{p}}\hat{t}^{\frac{1}{p}}\frac{\hat{\pi}}{2})}{\Si((-1)^{\frac{1}{p}}\hat{t}^{\frac{1}{p}}\frac{\hat{\pi}}{2})}
=\hat{t}^{\frac{1}{p}}\frac{\Sh'(\hat{t}^{\frac{1}{p}}\frac{\hat{\pi}}{2})}{\Sh(\hat{t}^{\frac{1}{p}}\frac{\hat{\pi}}{2})}\ ,$$
we have
\begin{eqnarray*}
f'(t)&=&-\frac{1}{p}\hat{t}^{\frac{1-p}{p}}\frac{\Sh'(\hat{t}^{\frac{1}{p}}\frac{\hat{\pi}}{2})}{\Sh(\hat{t}^{\frac{1}{p}}\frac{\hat{\pi}}{2})}
+\hat{t}^{\frac{1}{p}}(-\frac{1}{p}\hat{t}^{\frac{1-p}{p}})\frac{\hat{\pi}}{2}\frac{\Sh''(\hat{t}^{\frac{1}{p}}\frac{\hat{\pi}}{2})\Sh(\hat{t}^{\frac{1}{p}}\frac{\hat{\pi}}{2})-\Sh'^{2}(\hat{t}^{\frac{1}{p}}\frac{\hat{\pi}}{2})   }{\Sh^{2}(\hat{t}^{\frac{1}{p}}\frac{\hat{\pi}}{2})}\ ,\\
&=&\frac{-\frac{1}{p}\hat{t}^{\frac{1-p}{p}}}{\Sh^{2}(\hat{t}^{\frac{1}{p}}\frac{\hat{\pi}}{2})}\left[\Sh'(\hat{t}^{\frac{1}{p}}\frac{\hat{\pi}}{2})\Sh(\hat{t}^{\frac{1}{p}}\frac{\hat{\pi}}{2})
+\frac{\hat{\pi}}{2}\hat{t}^{\frac{1}{p}}\left(\frac{\Sh^{p}(\hat{t}^{\frac{1}{p}}\frac{\hat{\pi}}{2})}{\Sh'^{p-2}(\hat{t}^{\frac{1}{p}}\frac{\hat{\pi}}{2})}-\Sh'^{2}(\hat{t}^{\frac{1}{p}}\frac{\hat{\pi}}{2})\right)
\right]\ ,\\
&=&\frac{-\frac{1}{p}\hat{t}^{\frac{1-p}{p}}}{\Sh^{2}(\hat{t}^{\frac{1}{p}}\frac{\hat{\pi}}{2})}  \left[
               \Sh'(\hat{t}^{\frac{1}{p}}\frac{\hat{\pi}}{2})\Sh(\hat{t}^{\frac{1}{p}}\frac{\hat{\pi}}{2})-\frac{\hat{\pi}}{2}\hat{t}^{\frac{1}{p}}\Sh'^{2-p}(\hat{t}^{\frac{1}{p}}\frac{\hat{\pi}}{2})\right]\\\
&\equiv&\frac{-\frac{1}{p}\hat{t}^{\frac{1-p}{p}}}{\Sh^{2}(\hat{t}^{\frac{1}{p}}\frac{\hat{\pi}}{2})}\tilde{g}(t)\ .
\end{eqnarray*}
Using similar argument as step 1, we can show $\tilde{g}(t)>0$ and hence $f'(t)<0$ for all $t<0$.

\item[3.] If $f(t)=-f(t-m)$, then $f'(t)\frac{dt}{dm}=-f'(t-m)(\frac{dt}{dm}-1)$ and $$\frac{dt}{dm}=\frac{f'(t-m)}{f'(t)+f'(t-m)}>0\ .$$
If we can show $f'(t_{2}-m)<f'(t_{2})$ and $f'(t_{1}-m)>f'(t_{1})$, then
$$\frac{dt_{1}}{dm}=\frac{f'(t_{1}-m)}{f'(t_{1})+f'(t_{1}-m)}<\frac{f'(t_{2}-m)}{f'(t_{2})+f'(t_{2}-m)}=\frac{dt_{2}}{dm}\ .$$
Hence $\frac{d}{dm}(t_{2}-t_{1})(m)>0$ for all $m>0$. Furthermore
$$(t_{2}-t_{1})(m)>\lim_{m\to 0^{+}}(t_{2}-t_{1})(m)=2^{p}-1\ .$$

\item[4.] First, note that $t_{2}>m$ and $f(t_{2})>0$ for $m<3^{p}-1$.
Since $f(t_{2})=-f(t_{2}-m)$ and
\begin{eqnarray*}
f'(t)
&=&\frac{1}{pt}(f(t)-\frac{\hat{\pi}}{2}(1+|\T(t^{\frac{1}{p}}\frac{\hat{\pi}}{2})|^{p})|f(t)|^{2})\ \\
&=&\frac{1}{pt}f(t)-\frac{\hat{\pi}}{2pt}(1+\frac{t}{|f(t)|^{p}})|f(t)|^{2}\ ,
\end{eqnarray*}
we have
\begin{eqnarray*}
f'(t_{2}-m)-f'(t_{2})&=&-\frac{f(t_{2})}{p}(\frac{1}{t_{2}-m}+\frac{1}{t_{2}})+\frac{\hat{\pi}|f(t_{2})|^{2}}{2p}(\frac{1}{t_{2}}(1+\frac{t_{2}}{|f(t_{2})|^{p}})-\frac{1}{t_{2}-m}(1+\frac{t_{2}-m}{|f(t_{2})|^{p}}))\ \\
&=&-\frac{(2t_{2}-m)f(t_{2})}{pt_{2}(t_{2}-m)}-\frac{m\hat{\pi}|f(t_{2})|^{2}}{2pt_{2}(t_{2}-m)}\ \\
&<&0\ .
\end{eqnarray*}

\item[5.] Note $f(t_{1})<0$ for $m>0$, and $t_{1}-m>1-m>0$ if $m<1$. Since
$$
f'(t)
=\frac{1}{pt}f(t)-\frac{\hat{\pi}}{2pt}(1+\frac{t}{|f(t)|^{p}})|f(t)|^{2}\ ,
$$
we have
\begin{eqnarray*}
f'(t_{1})&=&\frac{f(t_{1})}{pt_{1}}-\frac{\hat{\pi}|f(t_{1})|^{2}}{2pt_{1}}-\frac{\hat{\pi}|f(t_{1})|^{2-p}}{2p}\ ,\\
f'(t_{1}-m)&=&-\frac{f(t_{1})}{p(t_{1}-m)}-\frac{\hat{\pi}|f(t_{1})|^{2}}{2p(t_{1}-m)}-\frac{\hat{\pi}|f(t_{1})|^{2-p}}{2p}\ ,
\end{eqnarray*}
and hence
\begin{eqnarray}
\nonumber f'(t_{1})+\frac{\hat{\pi}|f(t_{1})|^{2-p}}{2p}   &=&\frac{f(t_{1})}{pt_{1}}-\frac{\hat{\pi}|f(t_{1})|^{2}}{2pt_{1}}\ <\ 0\ ,\\
f'(t_{1}-m)+ \frac{\hat{\pi}|f(t_{1})|^{2-p}}{2p}   &=&-\frac{f(t_{1})}{p(t_{1}-m)}-\frac{\hat{\pi}|f(t_{1})|^{2}}{2p(t_{1}-m)}\ =\ -\frac{f(t_{1})}{p(t_{1}-m)}(1+\frac{\hat{\pi}}{2}f(t_{1})).\
\label{eq2.6}
\end{eqnarray}
Since LHS in (\ref{eq2.6}) is finite for $0<m<3^{p}-1$, $t_{1}$ is increasing in $m$, and $f(t)$ is decreasing in $t$, there exists unique $m^{*}$ such that $$t_{1}(m^{*})-m^{*}=0,\ 1+\frac{\hat{\pi}}{2}f(t_{1}(m^{*}))=0\ .$$ Hence,
$$
\begin{array}{lll}
t_{1}(m)-m
>0\ ,&1+\frac{\hat{\pi}}{2}f(t_{1}(m))>0&\ \mbox{on}\ (0,m^{*})\ ,\\
 t_{1}(m)-m
<0\ ,&1+\frac{\hat{\pi}}{2}f(t_{1}(m))<0&\ \mbox{on}\ (m^{*},3^{p}-1)\ .\\
\end{array}
$$
Furthermore,
$$f'(t_{1}-m)+ \frac{\hat{\pi}|f(t_{1})|^{2-p}}{2p} =-\frac{f(t_{1})}{p(t_{1}-m)}-\frac{\hat{\pi}|f(t_{1})|^{2}}{2p(t_{1}-m)}>0.\ $$
This implies
$$f'(t_{1}-m)>- \frac{\hat{\pi}|f(t_{1})|^{2-p}}{2p}> f'(t_{1})\ .$$
\end{enumerate}
\end{proof}

\newtheorem{th2.4}[th2.2]{Lemma}
\begin{th2.4}
\label{th2.4}

Let $s_{1}$ and $s_{2}$ be the first two zeros of $m\tan_{p}s=-\tan_{p}(sm)$ for $m>1$. Then
\begin{equation}
\frac{s_{2}(m)}{s_{1}(m)}>2\ .\label{eq2.2}
\end{equation}

\end{th2.4}

\begin{proof}

To do this, we claim that $$\frac{d}{dm}\frac{s_{2}(m)}{s_{1}(m)}=\frac{s'_{2}(m)s_{1}(m)-s_{2}(m)s'_{1}(m)}{s_{1}^{2}(m)}>0\ .$$
We first observe that, if $m\tan_{p}s=-\tan_{p}(sm)$, then $$\tan_{p}s+m\left(1+|\tan_{p}s|^{p}\right)\frac{ds}{dm}=-\left(1+|\tan_{p}(sm)|^{p}\right)\left(s+m\frac{ds}{dm}\right)\ ,$$
or equivalently
$$\frac{ds}{dm}=-\frac{\tan_{p}s+s(1+|\tan_{p}(sm)|^{p})}{m(1+|\tan_{p}s|^{p}+1+|\tan_{p}(sm)|^{p})}\ .$$
Hence$$s'_{2}(m)s_{1}(m)-s_{2}(m)s'_{1}(m)=\frac{F(m,s_{1},s_{2})}{m\tilde{F}(s_{1})\tilde{F}(s_{2})}\ ,$$
where
$$\tilde{F}(s)=1+|\tan_{p}s|^{p}+1+|\tan_{p}(sm)|^{p}\ ,$$
and
\begin{eqnarray*}
F(m,s_{1},s_{2})&=&s_{2}(\tan_{p}s_{1}+s_{1}(1+|\tan_{p}(s_{1}m)|^{p}))(1+|\tan_{p}s_{2}|^{p}+1+|\tan_{p}(s_{2}m)|^{p})\\
&&\ -s_{1}(\tan_{p}s_{2}+s_{2}(1+|\tan_{p}(s_{2}m)|^{p}))(1+|\tan_{p}s_{1}|^{p}+1+|\tan_{p}(s_{1}m)|^{p})\ \\
&=&(\tan_{p}s_{1}-s_{1}(1+|\tan_{p}s_{1}|^{p}))(\tan_{p}s_{2}+s_{2}(1+|\tan_{p}(s_{2}m)|^{p}))\\
&&\ -(\tan_{p}s_{2}-s_{2}(1+|\tan_{p}s_{2}|^{p}))(\tan_{p}s_{1}+s_{1}(1+|\tan_{p}(s_{1}m)|^{p}))\ \\
&=&(\tan_{p}s_{1}-s_{1}(1+|\tan_{p}s_{1}|^{p}))(\tan_{p}s_{2}+s_{2}(1+m^{p}|\tan_{p}s_{2}|^{p}))\\
&&\ -(\tan_{p}s_{2}-s_{2}(1+|\tan_{p}s_{2}|^{p}))(\tan_{p}s_{1}+s_{1}(1+m^{p}|\tan_{p}s_{1}|^{p}))\ .
\end{eqnarray*}
The last equality is because $m\tan_{p}s_{i}=-\tan_{p}(s_{i}m), i=1,2$.
Define
\begin{eqnarray*}
g_{1}(s)&=&\tan_{p}s-s(1+|\tan_{p}s|^{p})\ ,\\
g_{2}(s)&=&\tan_{p}s+s(1+m^{p}|\tan_{p}s|^{p})\ .
\end{eqnarray*}
Note that $g_{2}(s)>0$ for $s\in (0,\hat{\pi})$. Denote by $G(s)=\frac{g_{1}(s)}{g_{2}(s)}$. Since $$
\begin{array}{l}
\lim_{s\to 0}G(s)=  0\ ,\quad
\lim_{s\to \hat{\pi}}G(s)=-1\ ,\\
\lim_{s\to \frac{\hat{\pi}}{2}^{+}}G(s)=\lim_{s\to \frac{\hat{\pi}}{2}^{-}}G(s)=-\frac{1}{m^{p}}\ ,
\end{array}
$$
the function $G(s)$ is well-defined on $[0,\hat{\pi}]$.
\begin{enumerate}
  \item[1.]
For $m>3$, we have $s_{1}, s_{2}\in (0,\frac{\hat{\pi}}{2})$. If we can show  $G(s)$ is decreasing on $(0,\frac{\hat{\pi}}{2})$, then $\frac{g_{1}(s_{1})}{g_{2}(s_{1})}>\frac{g_{1}(s_{2})}{g_{2}(s_{2})}$ and hence $$\frac{d}{dm}\frac{s_{2}(m)}{s_{1}(m)}>0 \ .$$
Since, when $m\to 3^{+}$, $t_{2}(m)\to\frac{\hat{\pi}}{2}^{-}$ and $t_{1}\in (\frac{\hat{\pi}}{6},\frac{\hat{\pi}}{4})$, we have, for $m>3$,
$$\frac{s_{2}(m)}{s_{1}(m)}>\lim_{m\to 3^{+}}\frac{s_{2}(m)}{s_{1}(m)}>\frac{\hat{\pi}/2}{\hat{\pi}/4}=2\ .$$

Now, for $s\in (0,\frac{\hat{\pi}}{2})$, we have $\tan_{p}s > 0$ and
\begin{eqnarray*}
&&g'_{1}(s)g_{2}(s)-g_{1}(s)g'_{2}(s)\\
&&\ =m^{p}|\tan_{p}s|^{p-1}\left[s(1+|\tan_{p}s|^{p})(ps+(1-p)\tan_{p}s)-|\tan_{p}s|^{2}\right]\\
&&\quad +(1+|\tan_{p}s|^{p})\left[s(2+(1-p)|\tan_{p}s|^{p}-ps|\tan_{p}s|^{p-1})-\tan_{p}s\right]-\tan_{p}s\\
&&\ \equiv m^{p}|\tan_{p}s|^{p-1}G_{1}(s)+G_{2}(s)\ .
\end{eqnarray*}
\begin{enumerate}
  \item[(1)] Since $G_{2}(0)=0$ and, for $s\in (0,\frac{\hat{\pi}}{2})$,
\begin{eqnarray*}
G'_{2}(s)&=&p(1+|\tan_{p}s|^{p})|\tan_{p}s|^{p-2}[|\tan_{p}s|^{p}s((1-p)\tan_{p}s-ps)\\
 &&\ -p|\tan_{p}s|^{2} +(1-p)s(1+|\tan_{p}s|^{p})(\tan_{p}s+s) ]\ \\
&<&0\ ,
\end{eqnarray*}
we have $G_{2}(s)<0$ for $s\in (0,\frac{\hat{\pi}}{2})$.
  \item[(2)] For $s\in (0,\frac{\hat{\pi}}{2})$, we have $G_{1}(0)=0$,
\begin{eqnarray*}
G'_{1}(s)&=&(1+|\tan_{p}s|^{p})\left[ (1+p)(s-\tan_{p}s) +s|\tan_{p}s|^{p-1}(p^{2}s+(1-p^{2})\tan_{p}s)\right]\\\
&\equiv&(1+|\tan_{p}s|^{p})\hat{G}_{1}(s)\ ,
\end{eqnarray*}
and $\hat{G}_{1}(0)=0$,
\begin{eqnarray*}
\hat{G}'_{1}(s)&=& |\tan_{p}s|^{p-2}[ -p(1+p)|\tan_{p}s|^{2}+ps\tan_{p}s(1+2p-p^{2})\\
&&\ +(p-1)s|\tan_{p}s|^{p}(p^{2}s-p(1+p)\tan_{p}s)+(p-1)p^{2}s^{2}  ]\ \\
&\leq& |\tan_{p}s|^{p-2}[ -p(1+p)s\tan_{p}s+ps\tan_{p}s(1+2p-p^{2})\\
&&\ +(p-1)s|\tan_{p}s|^{p}(p^{2}s-p(1+p)\tan_{p}s)+(p-1)p^{2}s^{2}  ]\\\
&=&p(p-1)s|\tan_{p}s|^{p-2}\left[ p(s-\tan_{p}s)+p|\tan_{p}s|^{p}(s-\tan_{p}s)-|\tan_{p}s|^{p+1} \right]\ \\
&\leq&0\ .
\end{eqnarray*}
Hence $G_{1}(s)<0$ on $(0,\frac{\hat{\pi}}{2})$.
\end{enumerate}
This implies $g'_{1}(s)g_{2}(s)-g_{1}(s)g'_{2}(s)<0$ on $(0,\frac{\hat{\pi}}{2})$. Furthermore, $$G'(s)=\frac{g'_{1}(s)g_{2}(s)-g_{1}(s)g'_{2}(s)}{g_{1}^{2}(s)}<0\ .$$
That is $G(s)$ is decreasing on $(0,\frac{\hat{\pi}}{2})$.

\item[2.] For $m<3$,  we have $s_{2}\in (\frac{\hat{\pi}}{2},\hat{\pi})$, $\tan_{p}s_{2} <0 $ and $$
\frac{g_{1}(s_{2})}{g_{2}(s_{2})}=\frac{\tan_{p}s_{2}-s_{2}(1+|\tan_{p}s_{2}|^{p})}{\tan_{p}s_{2}+s_{2}(1+m^{p}|\tan_{p}s_{2}|^{p})}< -\frac{1}{m^{p}}\ ,
$$
since
\begin{eqnarray*}
0&>&m^{p}(\tan_{p}s_{2}-s_{2}(1+|\tan_{p}s_{2}|^{p}))+tan_{p}s_{2}+s_{2}(1+m^{p}|\tan_{p}s_{2}|^{p})\\\
&=&(m^{p}+1)\tan_{p}s_{2} +(1-m^{p})s_{2}\ ,
\end{eqnarray*}
is a tautology. Hence
$$\frac{g_{1}(s_{2})}{g_{2}(s_{2})}< -\frac{1}{m^{p}}=\lim_{s\to\frac{\hat{\pi}}{2}^{-}}G(s)<\frac{g_{1}(s_{1})}{g_{2}(s_{1})}\ ,$$
or equivalently $F(m,s_{1},s_{2})>0$. This implies $\frac{d}{dm}\frac{s_{2}(m)}{s_{1}(m)}>0$ for $m<3$. Furthermore
$$\frac{s_{2}(m)}{s_{1}(m)}>\lim_{m\to 1^{+}}\frac{s_{2}(m)}{s_{1}(m)} = \frac{\hat{\pi}}{\hat{\pi}/2}=2\ .$$
\end{enumerate}
\end{proof}

\section{Proof of Main Theorem}

\begin{proof}[Proof of Theorem  \ref{th1.1}]

For $M>0$, denote $$A_{M}=\{ 0\leq q(x)\leq M : q\
\mbox{is single-well with a transition point at}\ \frac{\hat{\pi}}{2}\}.$$ Let
$E[q]= (\lambda_2-\lambda_1)[q]$. Then $E[q]$ is bounded on $A_{M}$ and, hence, $E[q]$ attains its minimum  at some $q_{0}$ in $A_{M}$.
For $q(x)\in A_{M}$, define by $q(x,t)=tq(x)+(1-t)q_{0}(x)$ the one-parameter family of
potentials , where $0<t<1$.

By Lemma \ref{th2.2}, there exist $0\leq x_{-}<x_{0}<x_{+}\leq \hat{\pi}$, such that $y_{2}(x_{0},0)=0$ and
$$|y_{2}(x,0)|^{p}-|y_{1}(x,0)|^{p}\left\{\begin{array}{l}>0\ \mbox{on}\ (0,x_{-})\cup (x_{+},\hat{\pi})\ ,\\
<0\ \mbox{on}\ (x_{-}, x_{+})\ .\end{array}\right.$$

\begin{enumerate}
\item[1.] Assume $x_{-}\leq \frac{\hat{\pi}}{2}< x_{+}$. Let $$q(x)=\left\{\begin{array}{l}q_{0}(x_{-})\ \mbox{on}\ (0,\frac{\hat{\pi}}{2})\ ,\\
q_{0}(x_{+})\ \mbox{on}\ (\frac{\hat{\pi}}{2},\hat{\pi})\ . \end{array}\right.$$
By the optimality of $q_{0}$, we have, using Lemma \ref{th2.1},
$$0\leq \frac{d}{dt}(\la_{2}(t)-\la_{1}(t))=\int_{0}^{\hat{\pi}}(q(x)-q_{0}(x))(|y_{2}(x,0)|^{p}-|y_{1}(x,0)|^{p})dx\leq 0 \ .$$
This implies $q_{0}=q(x)$.

\item[2.] Assume $\frac{\hat{\pi}}{2}<x_{-}$ (the case for $x_{+}<\frac{\hat{\pi}}{2}$ is similar). Let
$$q(x)=\left\{\begin{array}{l}0\ \mbox{on}\ (0,x_{-})\ ,\\
M\ \mbox{on}\ (x_{-},\hat{\pi})\ . \end{array}\right.$$
Since $y_{n}(x,0)$ is normalized, we have
\begin{eqnarray*}
\int_{0}^{x_{-}} (|y_{2}(x,0)|^{p}-|y_{1}(x,0)|^{p})dx &>&0\ ,\\
\int_{x_{-}}^{\hat{\pi}} (|y_{2}(x,0)|^{p}-|y_{1}(x,0)|^{p})dx &<&0\ .
\end{eqnarray*}

By the optimality of $q_{0}$, we have
\begin{eqnarray*}
0&\leq& \frac{d}{dt}(\la_{2}(t)-\la_{1}(t))=\int_{0}^{\hat{\pi}}(q(x)-q_{0}(x))(|y_{2}(x,0)|^{p}-|y_{1}(x,0)|^{p})dx\\
&=&-q_{0}(\frac{\hat{\pi}}{2})\int_{0}^{x_{-}} (|y_{2}(x,0)|^{p}-|y_{1}(x,0)|^{p})dx+(M-q_{0}(x_{+}))\int_{x_{-}}^{\hat{\pi}} (|y_{2}(x,0)|^{p}-|y_{1}(x,0)|^{p})dx\\\
&\leq&0\ .
\end{eqnarray*}
The only possibility is $q_{0}=q$. But in this case, the second eigenfunction can be expressed by
$$y_{2}(x)=\left\{\begin{array}{l}c\Si(\la_{2}^{\frac{1}{p}}x)\ \mbox{on}\ (0,\frac{\hat{\pi}}{2})\ ,\\
d\Si((\la_{2}-M)^{\frac{1}{p}}(\hat{\pi}-x))\ \mbox{on}\ (\frac{\hat{\pi}}{2},\hat{\pi})\ . \end{array}\right.$$
Since $\frac{\hat{\pi}}{2}<x_{-}<x_{0}<x_{+}$, we have $\la_{2}^{\frac{1}{p}}\frac{\hat{\pi}}{2} < \hat{\pi}$ and $(\la_{2}-M)^{\frac{1}{p}}\frac{\hat{\pi}}{2}>\hat{\pi}$. Furthermore,
$$(\la_{2}-M)^{\frac{1}{p}} > \la_{2}^{\frac{1}{p}}\ .$$
This is impossible and hence this case is refused.
\end{enumerate}

By above discussion, we may assume  $$q_{0}(x)=\left\{\begin{array}{l}m\ \mbox{on}\ (0,\frac{\hat{\pi}}{2})\ ,\\
0\ \mbox{on}\ (\frac{\hat{\pi}}{2},\hat{\pi})\ . \end{array}\right.$$
In this case, the eigenfunction corresponding to the eigenvalue $\la$ can be expressed as
 $$y(x)=\left\{\begin{array}{l}c\Si(\la^{\frac{1}{p}}x)\ \mbox{on}\ (0,\frac{\hat{\pi}}{2})\ ,\\
d\Si((\la-m)^{\frac{1}{p}}(\hat{\pi}-x))\ \mbox{on}\
(\frac{\hat{\pi}}{2},\hat{\pi})\ . \end{array}\right.$$ Here,
$\la$ is an eigenvalue if $\la$ is a solution of
$$\frac{\la^{\frac{1}{p}}\Si'(\la^{\frac{1}{p}}\frac{\hat{\pi}}{2})}{\Si(\la^{\frac{1}{p}}\frac{\hat{\pi}}{2})}=-\frac{(\la-m)^{\frac{1}{p}}\Si'((\la-m)^{\frac{1}{p}}\frac{\hat{\pi}}{2})}{\Si((\la-m)^{\frac{1}{p}}\frac{\hat{\pi}}{2})}\
,$$ or equivalently $$\la^{\frac{1}{p}}\cot_{p}
(\la^{\frac{1}{p}}\frac{\hat{\pi}}{2})
=-(\la-m)^{\frac{1}{p}}\cot_{p}
((\la-m)^{\frac{1}{p}}\frac{\hat{\pi}}{2})\ .$$ By Lemma
\ref{th2.3}, we obtain the eigenvalue gap $\la_{2}-\la_{1}\geq
2^{p}-1$ and the equality holds if and only if $q$ is constant.

Finally, we assume
$$q(x,t)=\left\{\begin{array}{l}t\ \mbox{on}\ (0,a)\ ,\\
0\ \mbox{on}\ (a,\hat{\pi})\ , \end{array}\right.$$
for $t\geq 0$.
Then $y_{1}(x,0)=(\frac{p}{\hat{\pi}})^{\frac{1}{p}}\sin_{p}x, y_{2}(x,0)=(\frac{p}{\hat{\pi}})^{\frac{1}{p}}\sin_{p}(2x)$ and $\int_{0}^{\frac{\hat{\pi}}{2}}(|y_{2}(x,0)|^{p}-|y_{1}(x,0)|^{p})dx=0$. Hence $$\frac{d}{dt}(\la_{2}-\la_{1})(0)=\int_{0}^{a}(|y_{2}(x,0)|^{p}-|y_{1}(x,0)|^{p})dx<0\ ,$$
for $0<a- \frac{\hat{\pi}}{2}<<1$. Furthermore, for small $t>0$, we have $(\la_{2}-\la_{1})(t)<(\la_{2}-\la_{1})(0)=2^{p}-1$.
\end{proof}

\begin{proof}[Proof of Theorem  \ref{th1.2}.]
{\bf Part (a).} For $M>1$, denote $$A_{M}=\{ \frac{1}{M}\leq
\rho(x)\leq M : \rho\ \mbox{is single-barrier with a transition
point at}\ \frac{\hat{\pi}}{2}\}.$$ Let $R[q]=
\frac{\mu_{2}}{\mu_{1}}[q]$. Then $R[q]$ is bounded on $A_{M}$ and,
hence, $R[q]$ attains its minimum  at some $\rho_{0}$ in $A_{M}$.
For $\rho(x)\in A_{M}$, define $\rho(x,t)=t\rho(x)+(1-t)\rho_{0}(x)$
be the one-parameter family of densities, where $0<t<1$. Similar to
the proof of Theorem \ref{th1.1}, it can be showed that the optimal
$\rho_{0}$ must have the form
$$\rho_{0}=\left\{\begin{array}{l}1\ \mbox{on}\ (0,\frac{\hat{\pi}}{2})\ ,\\ L\ \mbox{on}\ (\frac{\hat{\pi}}{2},\hat{\pi}) \ ,\end{array}\right.\quad \mbox{or}\quad \rho_{0}=\left\{\begin{array}{l}L\ \mbox{on}\ (0,\frac{\hat{\pi}}{2})\ ,\\ 1\ \mbox{on}\ (\frac{\hat{\pi}}{2},\hat{\pi}) \ , \end{array}\right.$$ for some $L\geq 1$. W.L.O.G., we only discuss the first case. In this case, the eigenfunction corresponding to the eigenvalue $\mu$ can be expressed as
 $$y(x)=\left\{\begin{array}{l}c\Si(\mu^{\frac{1}{p}}x)\ \mbox{on}\ (0,\frac{\hat{\pi}}{2})\ ,\\
d\Si((\mu L)^{\frac{1}{p}}(\hat{\pi}-x))\ \mbox{on}\ (\frac{\hat{\pi}}{2},\hat{\pi})\ . \end{array}\right.$$
Here, $\mu$ is an eigenvalue if $\mu$ is a solution of
$$\frac{\mu^{\frac{1}{p}}\Si'(\mu^{\frac{1}{p}}\frac{\hat{\pi}}{2})}{\Si(\mu^{\frac{1}{p}}\frac{\hat{\pi}}{2})}=-\frac{(\mu L)^{\frac{1}{p}}\Si'((\mu L)^{\frac{1}{p}}\frac{\hat{\pi}}{2})}{\Si((\mu L)^{\frac{1}{p}}\frac{\hat{\pi}}{2})}\ ,$$
or equivalently
$$\mu^{\frac{1}{p}}\cot_{p}(\mu^{\frac{1}{p}}\frac{\hat{\pi}}{2})=-(\mu L)^{\frac{1}{p}}\cot_{p}((\mu L)^{\frac{1}{p}}\frac{\hat{\pi}}{2})\ .$$
Let $m=L^{\frac{1}{p}}$ and $s=\mu^{\frac{1}{p}}\frac{\hat{\pi}}{2}$. Then we obtain $$m\tan_{p}s=-\tan_{p}(sm)\ .$$
By  Lemma \ref{th2.4}, we obtain the eigenvalue ratio $\frac{\mu_{2}}{\mu_{1}}\geq 2^{p}$ and the equality holds if and only if $\rho$ is constant.

Finally, we assume
$$\rho(x,t)=\left\{\begin{array}{l}t\ \mbox{on}\ (0,a)\ ,\\
1\ \mbox{on}\ (a,\hat{\pi})\ , \end{array}\right.$$
for $t\geq 1$.
Then $y_{1}(x,1)=(\frac{p}{\hat{\pi}})^{\frac{1}{p}}\sin_{p}x, y_{2}(x,1)=(\frac{p}{\hat{\pi}})^{\frac{1}{p}}\sin_{p}(2x)$ and $\int_{0}^{\frac{\hat{\pi}}{2}}(|y_{1}(x,1)|^{p}-|y_{2}(x,1)|^{p})dx=0$. Hence $$\frac{d}{dt}(\frac{\mu_{2}}{\mu_{1}})(1)=\frac{\mu_{2}(1)}{\mu_{1}(1)}\int_{0}^{a}(|y_{1}(x,1)|^{p}-|y_{2}(x,1)|^{p})dx<0\ ,$$
for $0< \frac{\hat{\pi}}{2}-a<<1$. Furthermore, for small $t>0$, we have $(\frac{\mu_{2}}{\mu_{1}})(t)<(\frac{\mu_{2}}{\mu_{1}})(1)=2^{p}$.

\noindent {\bf Part (b).} We give an alternative proof with respect
to part (a). Consider the one-parameter family of densities
$\rho(x,t)=t\rho(x)+(1-t)\epsilon$, where $0<t<1$ and $\epsilon$ is
a positive constant.
Denote by $\{\mu_n(t),\ y_n(x,t)\}$ the $n$-th normalized eigenpair corresponding to the density
$\rho(x,t)$.
By Lemma \ref{th2.2}, there are points $x_\pm(t)$
with
$$
0<x_-(t)<\frac{\hat{\pi}}{2}<x_+(t)<\hat{\pi},~~~~x_-(t)+x_+(t)=\hat{\pi}
$$
such that
\begin{equation}
\begin{cases}|y_2(x,t)|^p>|y_1(x,t)|^p~~~\hbox{on}~~(0,x_-(t))\cup (x_+(t),\hat{\pi}),\\
|y_2(x,t)|^p<|y_1(x,t)|^p~~~\hbox{on}~~(x_-(t),x_+(t)).\\
\end{cases}
\label{eq3.2}\\
\end{equation}
Now, we claim that $$\frac{d}{dt}[\frac{\mu_2(t)}{\mu_1(t)}]\le
0~~~\mbox{for}~~~0<t<1.$$

From Lemma \ref{th2.2}, we have
\begin{equation*}
\frac{d}{dt}[\frac{\mu_2(t)}{\mu_1(t)}]=\frac{\mu_2(t)}{\mu_1(t)}\int_0^{\hat{\pi}}(\rho(x)-\epsilon)
[|y_1(x,t)|^p-|y_2(x,t)|^p]dx.
\end{equation*}
Since $\rho(x)$ is a symmetric single-well density, we obtain
\begin{eqnarray*}
\int_0^{\hat{\pi}}\rho(x)[|y_1(x,t)|^p-|y_2(x,t)|^p]dx&=&\int_{(0,x_-(t))\cup
(x_+(t),a)}\rho(x) [|y_1(x,t)|^p-|y_2(x,t)|^p]dx\nonumber\\
&+&\int_{x_-(t)}^{x_+(t)}\rho(x)
[|y_1(x,t)|^p-|y_2(x,t)|^p]dx\nonumber\\
&\le&\rho(x_-(t))\int_0^{\hat{\pi}}[|y_1(x,t)|^p-|y_2(x,t)|^p]dx.
\end{eqnarray*}
So,
\begin{equation}
\int_0^{\hat{\pi}}(\rho(x)-\epsilon)[|y_1(x,t)|^p-|y_2(x,t)|^p]dx\le[\rho(x_-(t))-\epsilon)]\int_0^{\hat{\pi}}[|y_1(x,t)|^p-|y_2(x,t)|^p]dx.
\label{eq3.3}
\end{equation}
The normalization condition
$\int_0^{\hat{\pi}}[t\rho(x)+(1-t)\epsilon]|y_n(x,t)|^pdx=1$ gives
\begin{equation}
\int_0^{\hat{\pi}}(|y_1(x,t)|^p-|y_2(x,t)|^p)dx=\frac{t}{\epsilon}\int_0^{\hat{\pi}}[\rho(x)-\epsilon](|y_2(x,t)|^p-|y_1(x,t)|^p)dx.
\label{eq3.4}
\end{equation}
So, by (\ref{eq3.3}), we obtain
\begin{equation*}
[\frac{\rho(x_-(t))t}{\epsilon}+(1-t)]\int_0^{\hat{\pi}}[\rho(x)-\epsilon](|y_1(x,t)|^p-|y_2(x,t)|^p)dx\le
0.
\end{equation*}
Since $0<t<1$, this implies that
\begin{equation}
\int_0^{\hat{\pi}}[\rho(x)-\epsilon](|y_1(x,t)|^p-|y_2(x,t)|^p)dx\le 0,\label{eq3.5}
\end{equation}
from which it follows that
$$\frac{d}{dt}[\frac{\mu_2(t)}{\mu_1(t)}]\le 0~~~for~~~0<t<1.$$

Finally, by the continuity of eigenvalues, we obtain
\begin{equation*}
\frac{\mu_2[\rho]}{\mu_1[\rho]}=\frac{\mu_2(1)}{\mu_1(1)}\le\frac{\mu_2(0)}{\mu_1(0)}
=\frac{\mu_2[\epsilon]}{\mu_1[\epsilon]}=2^p.
\end{equation*}
The equality occurs only if $\frac{\mu_2(t)}{\mu_1(t)}$
is a constant. In this case, the equality holds in (\ref{eq3.5}),
and it follows from (\ref{eq3.4}) that
\begin{equation*}
\int_0^{\hat{\pi}}\rho(x)(|y_1(x,t)|^p-|y_2(x,t)|^p)dx=\int_0^{\hat{\pi}}(|y_1(x,t)|^p-|y_2(x,t)|^p)dx=0.
\end{equation*}
This together with (\ref{eq3.3}) implies that $\rho(x)$ is a
constant a.e..
\end{proof}

\noindent{\bf Remark:}
\begin{enumerate}
\item[(i)] In Theorem 1.1, if we replace 'single-well' by
'single-barrier', the method fails because the inequality in Lemma
2.3 is the same. Thus the case for `single-barrier` potential is
still unknown.
\item[(ii)] In Theorem 1.2(a), if the condition 'single-barrier' is replaced by
'single-well', our proof can not work because the inequality in Lemma
2.4 remains the same. Thus the case for `single-well` densities of
p-Laplacian is still open.
\item[(iii)] In Theorem 1.2(b), if the condition 'symmetric single-well' is replaced by
'symmetric single-barrier', then the equality is reversed.
\end{enumerate}

\section{Appendix}

\begin{proof}[Proof of Lemma \ref{th2.1}]
In the following computation, we drop the suffix for convenience. Denote $\dot{y}=\frac{\partial y}{\partial t}$.
Differentiating (\ref{eq1.1}) with respect to t, we have
\begin{eqnarray*}
&~& (p-2)y'(x,t)^{(p-3)}\dot{y}'(x,t)y''(x,t)+ |y'(x,t)|^{p-2}\dot{y}''(x,t)\\
&+& \left(\dot{\lambda}(t)\rho(x,t)
+ \lambda(t)\dot{\rho}(x,t) - \dot{q}(x,t)\right)y(x,t)^{(p-1)}\\
&+&(p-1) (\lambda(t)\rho(x,t)-q(x))|y(x,t)|^{p-2}\dot{y}(x,t)=0\
 .
\end{eqnarray*}
Multiplying it by $y(x.t)$ and by (\ref{eq1.1}), we have that
\begin{eqnarray}
&&\left(-\dot{\lambda}(t)\rho(x,t) -\lambda(t)\dot{\rho}(x,t) +\dot{q}(x,t)\right)|y(x,t)|^p\nonumber\\
&&\ =\left((p-2)y'(x,t)^{(p-3)}y''(x,t)\dot{y}'(x,t)+|y'(x,t)|^{p-2}\dot{y}''(x,t)\right)y(x,t)-(y'(x,t)^{(p-1)})'\dot{y}(x,t)\ ,\nonumber\\
&&\ =\left(|y'(x,t)|^{p-2}\dot{y}'(x,t)\right)'y(x,t)-(y'(x,t)^{(p-1)})'\dot{y}(x,t)\ ,\nonumber\\
&&\ \equiv I - II\ .\label{eq.A1}
\end{eqnarray}
Since
\begin{eqnarray*}
\int_{0}^{\hat{\pi}}I&=&|y'(x,t)|^{p-2}\dot{y}'(x,t)]y(x,t)|_0^{\hat{\pi}}-\int_0^{\hat{\pi}}|y'(x,t)|^{p-2}\dot{y}'(x,t)y'(x,t)dx\ ,\\
&=&-\int_0^{\hat{\pi}}|y'(x,t)|^{p-2}\dot{y}'(x,t)y'(x,t)dx\ ,
\end{eqnarray*}
and
\begin{eqnarray*}
\int_{0}^{\hat{\pi}}II&=&y'(x,t)^{(p-1)}\dot{y}(x,t)|_0^{\hat{\pi}}-\int_0^{\hat{\pi}}|y'(x,t)|^{p-2}\dot{y}'(x,t)]y'(x,t)dx\ ,\\
&=&-\int_0^{\hat{\pi}}|y'(x,t)|^{p-2}\dot{y}'(x,t)y'(x,t)dx\ ,
\end{eqnarray*}
after integrating (\ref{eq.A1}) over $[0,\hat{\pi}]$ with respect
to x, it follows from $\int_0^{\hat{\pi}}\rho|y|^pdx=1$ that
$$\dot{\lambda}(t)=-\int_0^{\hat{\pi}}\lambda(t)\dot{\rho}(x,t)|y(x,t)|^pdx+\int_0^{\hat{\pi}}\dot{q}(x,t)|y(x,t)|^pdx\ .
$$
Let $\la=\la_{n}$. The proof is complete.

\end{proof}

\end{document}